\documentclass[11pt]{article}
\usepackage{amssymb, amsmath, amsthm, graphicx}
\usepackage{fullpage}
\usepackage{enumerate}

\newtheorem{theorem}{Theorem}[section]
\newtheorem{lemma}[theorem]{Lemma}
\newtheorem{corollary}[theorem]{Corollary}
\newtheorem{conjecture}[theorem]{Conjecture}

\author{Radoslav Fulek\thanks{The author gratefully acknowledge support from the Swiss National
     Science Foundation Grant No.  200021-125287/1 and ESF Eurogiga project GraDR as GA\v{C}R GIG/11/E023.}
\thanks{Department of Applied Mathematics, Charles University,  Prague, Email: {\tt radoslav.fulek@gmail.com}}}

\title{Estimating the number of disjoint edges in simple topological graphs via cylindrical drawings\footnote{The most of the content of the paper will appear in the proceedings of
Symposium on Computational Geometry 2013 as a part of the contribution entitled ``Topological Graphs: Empty Triangles and Disjoint Matchings''
co-authored by Radoslav Fulek and Andres J. Ruiz-Vargas.}}

\begin{document}
\maketitle
%
%

\begin{abstract}
A topological graph drawn on a cylinder whose base is horizontal  is \emph{angularly monotone}
if every vertical line intersects every edge at most once.
 Let $c(n)$ denote the maximum number $c$ such that every simple angularly monotone drawing of a complete graph on $n$ vertices contains at least
 $c$ pairwise disjoint edges.
We show that for every simple complete topological graph $G$ there exists $\Delta$, $0<\Delta<n$, such that $G$ contains at least $\max \{\frac n\Delta, c(\Delta)\}$ pairwise disjoint edges.
By combining our result with a result of T\'oth we obtain an alternative proof for the best known lower bound of $\Omega(n^\frac 13)$ on the maximum number
of pairwise disjoint edges in a simple complete topological graph proved by Suk.
Our proof is based on a result of Ruiz-Vargas.
\end{abstract}

%
%
%
%



\section{Introduction}

A \emph{topological graph} is a graph drawn on a surface so that its vertices are represented by points and its edges are represented by Jordan arcs connecting the respective endpoints. In the present note we usually assume that a graph is drawn in the plane unless we state otherwise.  A \emph{topological graph} is \emph{simple} if every pair of its edges meet at most once either in a common vertex or at a proper crossing.
We use the words ``vertex'' and ``edge'' in both contexts, when referring to the elements of an abstract graph and also when referring to their planar counterparts.
A graph is \emph{complete} if there is an edge between every pair of vertices. Throughout this note $n$ denotes the number of vertices in a  graph.

The notion of a simple topological graph is a natural relaxation of a \emph{geometric graph} in which the edges are represented by  straight-line segments, and hence,
every pair of them is allowed to meet at most once. In geometric graphs we do not allow overlapping edges or edges passing through a vertex.

A \emph{matching} in a graph $G$ is a set of edges in $G$ no two of which share a vertex.
Two edges in a topological graph are \emph{disjoint} if their corresponding Jordan arcs do not intersect.
A \emph{disjoint matching} in a topological graph is the matching consisting of pairwise disjoint edges.
The aim of this note is to relate the lower bound on the size of maximal disjoint matching in simple topological graphs
to the lower bound of the same quantity in certain cylindrical drawings of  simple topological graphs.

A topological graph drawn on a cylinder whose base is horizontal  is \emph{angularly monotone}
if every vertical line intersects every edge at most once.
Angularly monotone topological graphs are closely related to $x$-monotone topological graphs, since they can be represented as $x$-monotone topological graphs, if we allow edges to consist of two connected components one of which ends in positive infinity and the other one in negative infinity.
Our main result is the following.

\begin{theorem}
\label{thm:reduction}
 Let $c(n)$ denote the maximum number $c$ such that every cylindrical angularly monotone drawing
 of a simple complete graph on $n$ vertices contains a disjoint matching of size at least
 $c$.
For every simple complete topological graph $G$ there exists $\Delta$, $0<\Delta<n$, such that $G$ contains a disjoint matching of size at least $\max \{\frac n\Delta, c(\Delta)\}$.
\end{theorem}

We will show as a simple consequence of a result of T\'oth~\cite{T00} that every simple angularly monotone complete  topological graph on $n$ vertices contains a
disjoint matching of size $\Omega (\sqrt{n})$.
Combining this with our result we obtain a new and  shorter proof of a recent result of Suk~\cite{S12}.

\begin{theorem}\label{thm:main}~\cite{S12}
A complete simple topological graph on $n$ vertices drawn in the plane contains a disjoint matching of size  $\Omega(n^\frac 13)$.
\end{theorem}

Whereas the proof of Theorem~\ref{thm:main} from~\cite{S12} is based on the existence of a perfect matching with a low stabbing number
for set systems with polynomially bounded dual shattered function~\cite{CW89}, our reduction to cylindrical drawings
is based on an intuitive geometric argument which yields, for a complete simple topological graph, a crossing free subgraph with
minimum degree of at least two.
Moreover, Theorem~\ref{thm:reduction} has a potential to improve the bound in Theorem~\ref{thm:main} to $\sqrt{n}$.

It is very easy to see that a complete $x$-monotone topological graph contains a disjoint matching of size $\lfloor \frac n2 \rfloor$.
We conjecture that a similar statement is true for simple angularly monotone  topological graphs.

\begin{conjecture}
\label{conj:main}
A simple angularly monotone complete  topological graph on $n$ vertices contains a disjoint matching of size $\Omega (n)$.
\end{conjecture}

As mentioned previously the best corresponding lower bound we can prove is $\Omega(\sqrt n)$.
Note that any improvement of this lower bound would translate into a better bound in   Theorem~\ref{thm:main}.

\paragraph{Algorithmic aspects}
The previous theorem can be also thought of as a lower bound on the size of a largest independent set in the intersection graph
of edges in a complete simple topological graph.
Besides the fact that computing the maximum number of (pairwise) disjoint elements in an arrangement of geometric objects is
an old problem in computational geometry, this type of  research
is also motivated by an abundance of applications e.g. in frequency assigment~\cite{ELR04}, computational cartography~\cite{APS98}, or VLSI design~\cite{HM85}.
Determining the size of a largest independent set is NP-hard already
for intersection graphs of sets of segments in the plane lying in two directions~\cite{KN90}, disks~\cite{FPT81} and rectangles~\cite{AI83}.
Hence, known efficient algorithms searching for a large independent set in intersection graphs of geometric objects
can only approximate the size of a largest one for an overwhelming majority of types of geometric objects.

We remark that similarly as the proof of Suk our proof of Theorem~\ref{thm:main} gives an efficient algorithm for finding a disjoint matching of size $\Omega(n^\frac 13)$
in a complete simple topological graph,
 Thus, we have an efficient  $n^{\frac 23}$-approximation algorithm
for the problem of finding the maximum disjoint matching
in a complete simple topological graph, which beats the factor of $n^{1-\epsilon}$ for any small $\epsilon>0$ in the inapproximability  result for the independence number in general graphs due to Zuckerman~\cite{Z07}.

\paragraph{Disjoint matchings in simple topological graphs} Estimating the lower bound on the maximum disjoint matching
in dense graphs, i.e. in graphs with at least $\gamma n^2$ edges for a fixed constant $\gamma>0$, was so far less successful than in the case of complete graphs.
We are not aware of any sub-linear upper bound for the problem. The chain of several improvements~\cite{FS09,PST03,PT05} leaded, up to now, only to a lower bound of $\Omega(\log^{1+\epsilon} n)$~(\cite{FS09}, Corollary 1.11) for some small $\epsilon>0$,
 and  our proof of Theorem~\ref{thm:main}, as well as the one in~\cite{S12}, does not seem to be extendable to the case of dense graphs.
 The lower bound of $\Omega(\log^{1+\epsilon} n)$ was improved by Suk~\cite{S12+} to $n^{\Omega(\frac{1}{\log \log n})}$ in the case of dense graphs whose edges are drawn as $t$-monotone curves i.e. each edge is represented by a curve with at most $t-1$ vertical tangent points. In the previous bound the hidden constant in
 $\Omega$-notation depends on $t$ and $\gamma$.

  In sparse graphs, the research is focused on estimating the upper bound on the number of edges in a simple topological graph that does not contain a disjoint matching of size $k$.
 Here, the best known upper bound of $O(n\log^{4k-8} n)$, for an arbitrary $k>1$, is due to Pach and T\'oth~\cite{PT05}, while the right order of magnitude is believed to be linear in $n$ for
 any fixed $k$.
 We remark that in all these problems it is essential that the topological graph is simple, since it was shown in~\cite{PT05} that a complete graph can be drawn in the plane
so that every pair of its edges intersect either once or twice.

\medskip

The paper is organized as follows. In Section 2, we present  tools that are used in Section 3 to derive Theorem~\ref{thm:reduction} and~\ref{thm:main}.
In Section 4, we finish with some concluding remarks.

\section{Preliminaries}

First we present the result of Ruiz-Vargas~\cite{RV13+} stated in Corollary~\ref{strong} which is the main tool for proving Theorem~\ref{thm:main}.

We say that a topological graph is \emph{plane} if it is free of edge crossings.
Let $G$ be a simple topological graph and $C$ be a cycle in $G$. We say that $C$ is \emph{plane} if no two edges of  $C$ cross.
 We use $C$ to denote both the cycle and the closed curve formed by its edges.  For a subgraph $H$ of $G$ and a vertex $v$ not in $H$, we say that $e$ is an edge from $v$ to $H$ if $e$ is incident to $v$ and to a vertex from $H$.

If $C$ is a simple closed curve, $\mathbb R^2\backslash C$ is partitioned uniquely, into two connected sets, one bounded and the other one unbounded. The latter  is referred to as the \emph{exterior} of $C$, while the former is called the \emph{interior} of $C$.

\begin{lemma}~\cite{RV13+}\label{mainlemma}
  Let $G$ be a simple topological graph, $C$ be a plane cycle of $G$ and $v$ be a vertex of $G$ in the interior of $C$ (resp. in the exterior of $C$). Suppose that for every $c\in V(C)$ we have $vc\in E(G)$. Then there exist at least two edges from $v$ to $V(C)$ that are contained in the interior of $C$ (resp. in the exterior of $C$).
\end{lemma}

 A \emph{face} of a plane graph is a connected
component of its complement.
We say that a vertex is incident to a face $F$ if it is contained in the closure of $F$, but not in $F$. For a face $F$ and a vertex $v$ in the interior of $F$, we say that an edge $e$ from $v$ to a vertex incident to $F$ is \emph{contained} in $F$, if the relative interior of $e$ is a subset of $F$.
Using Lemma~\ref{mainlemma} it is not hard to derive the following.
  \begin{corollary}~\cite{RV13+}\label{strong}
  Let $G$ be a simple topological graph and $H$ be a connected plane subgraph of $G$. Let $v$ be a vertex of $G$ that is not in $H$, and $F$ be the face of $H$ that contains $v$. Assume that for every vertex $w$ incident to $F$ we have $vw\in E(G)$. Then there exist two edges in $G$ from $v$ to $F$ that are contained in $F$.
  \end{corollary}
 A topological graph $G$ is $x$\emph{-monotone} if every vertical line intersects every edge of $G$ at most once.
  A topological graph $G$ is \emph{quasi} $x$\emph{-monotone} if every vertical line passing through a vertex intersects every edge of $G$ at most once (including its endpoints). We remind the reader that a topological graph drawn so that all the edges are represented as straight-line segments is called a \emph{geometric graph}.
  In order to derive Theorem~\ref{thm:main} from Theorem~\ref{thm:reduction} we use a result of G. T\'{o}th~\cite{T00} stating that a geometric graph without a
  disjoint matching of size $k$ has at most $O(k^2n)$ edges.
    \begin{theorem}~\cite{T00}\label{thm:main2}
  A geometric graph on $n$ vertices without  a disjoint matching of size $k$ has at most $O(k^2n)$ edges.
  \end{theorem}
   In fact, this result can be extended to simple quasi $x$-monotone topological graphs essentially without changing its proof.
   Indeed, the only properties of geometric graphs used in the proof of Theorem~\ref{thm:main2}
   are simplicity and $x$-monotonicity. Then by replacing each edge $e$ in a quasi $x$-monotone simple topological graph $G$
   by the polygonal line whose segments connect two consecutive intersections of $e$ with the vertical lines passing through the vertices of $G$,
   we obtain a simple $x$-monotone topological graph in which two edges cross if and only if they cross in $G$.
\begin{lemma}\label{lemma:main2}
A simple quasi $x$-monotone topological graph on $n$ vertices without a disjoint matching of size $k$  has at most $O(k^2n)$ edges.
  \end{lemma}
Finally, in the proof of Theorem~\ref{thm:main2} the large disjoint matching is found as a longest chain (Theorem 2 in~\cite{T00}) in a geometrically defined poset
that can be  constructed efficiently. This fact will easily imply that our proof of Theorem~\ref{thm:main} gives a polynomial time algorithm for the corresponding problem, since the algorithmic task of finding a longest chain in a poset can be carried out in a quadratic time (see e.g.~\cite{SPIN03}, Chapter 8) in the size
of the ground set.

\section{Disjoint matchings}

\subsection{Proof of Theorem~\ref{thm:reduction}}
\begin{figure}
\centering
\includegraphics[width=.8 \textwidth]{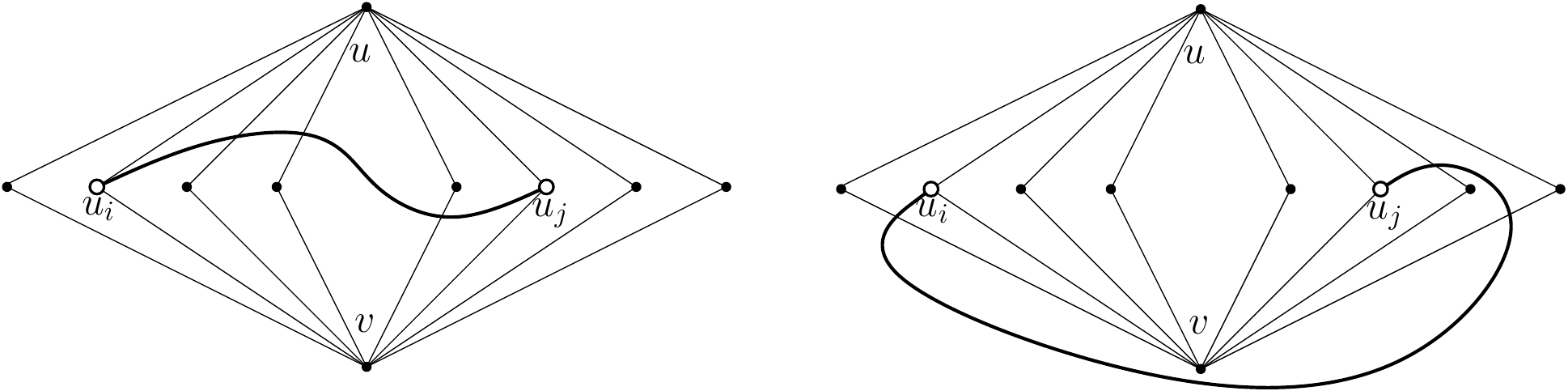}
\caption{Left: The edge $u_iu_j$ is drawn in the interior of the cycle $uu_ivu_j$. Right: The edge $u_iu_j$ is drawn in the exterior of the cycle $uu_ivu_j$.}
\label{diamonds}
\end{figure}

If $H$ is a  subgraph of a topological graph $G$ and $E'\subseteq E(G)$,
we denote by $H+E'$ the topological  subgraph of $G$ with the vertex set $V(H) \cup V(E')$ and the edge set $E(H) \cup E'$.

Let $v$ be a vertex of a complete simple topological graph $G$. Using Corollary~\ref{strong} we will show that there exists a connected plane subgraph $G'$
of $G$ containing all the edges incident to $v$ with minimum vertex degree two (and hence with at least $n-1+\lfloor \frac n2 \rfloor$ edges).

We inductively construct a sequence of plane subgraphs of $G$, $G_0\subset G_1 \subset \ldots \subset G_r=G'$ such that $G_0$ consists of all the edges incident
to $v$. Suppose that $G_i$ contains a vertex $u$ of degree one.  Let $e$ denote the unique edge incident to $u$ in $G_i$.
Let $F$ denote the unique face of $G_i-e$ such that $u$ is contained in $F$.
By applying Corollary~\ref{strong} on $G_i-e$  we obtain two edges $e'$ and $e''$ in $G$ incident to $u$ such that $G_i+\{e\} \cup \{e',e''\}$ is still a plane graph.  Indeed, the edges $e,e'$ and $e''$ are all incident to $u$, and thus, they do not cross each other as $G$ is simple.
 Since at least one of the edges $e'$ and $e''$ is different from $e$ we can set $G_{i+1}=G_i+\{e\} \cup \{e',e''\}$.

Let $\Delta$ denote the maximum degree of a vertex in $G'$ different from $v$.
In the rest of the proof we select a disjoint matching of $G$ in two ways. First, obtaining a matching of
size at least $\frac n\Delta$, and second, of size $c(\Delta)$. Clearly, once we establish this, the theorem is proved.

The first matching is obtained by taking a matching $M$ in $G'$ with the maximum
cardinality such that no edge in $M$ is incident to $v$. It is easy to see that $M$ has at least $\Omega(\frac n\Delta)$ edges. Indeed, the number of edges in $G'$ sharing a vertex with an edge in $M$ is
$O(|M|\Delta)$. Thus, if $|M|=o(\frac n\Delta)$ we would be able to find an edge in $G'$ that is not incident to $v$ and
that can be added to $M$ thereby obtaining a larger matching in $G'$ (for this recall that $G'$ has at least $\lfloor \frac n2 \rfloor$ edges not incident to $v$).

On the other hand, we can select a large disjoint matching as follows. Let $u$, $u\not=v$, denote a vertex in $G'$ of degree $\Delta$. Let $v_0,\ldots, v_{\Delta-1}$
 denote the neighbors of $u$ in $G'$ listed as the corresponding end-pieces of the edges $uv_i$ appear around $u$ in the clockwise order.
  By the fact that $G'$ is plane, this is also the counter-clockwise order in which the edges $vv_i$ appear around $v$.
   Let $P_0,P_1,\ldots,P_{\Delta-1}$ denote the paths $uv_0v,uv_1v,\ldots,  uv_{\Delta-1}v$, respectively.
Let $I_{i,j}$ and $I_{i,j}'$, $0\leq i\leq j< \Delta$, denote the cyclic intervals of integers $i+1, \ldots, j-1$ and $j+1,\ldots, \Delta-1,0,1,\ldots ,i-1$, respectively.
We suppose in the sequel that $i<j$. We say that an edge crosses a path if it crosses at least one of its edges.
 We claim the following (see Figure~\ref{diamonds} for an illustration).   \\

 (*) The edge $v_iv_j\in E(G)$ crosses  either all the paths $P_l, l\in I_{i,j}$, exactly once or all the paths $P_l, l\in I_{i,j}'$, exactly once. \\

 Indeed, since $G$ is simple, $v_iv_j$ crosses neither $P_i$ nor $P_j$, and hence, the relative interior of $v_iv_j$ is contained in a connected component of the complement of the Jordan arc corresponding  to the concatenation of $P_i$ and $P_j$. It follows that $v_iv_j$ crosses either all the paths $P_l, l\in I_{i,j}$, or all the paths $P_l, l\in I_{i,j}'$.
  As the two end-pieces of $v_iv_j$ lie in different connected components of the complement of the Jordan arc
 corresponding  to the concatenation of $P_i$ and $P_l$, $v_iv_j$ crosses $P_l$  an odd number of times. It follows that $v_iv_j$ crosses $P_l$ exactly once.

 Let $G_0'$ denote the subgraph of $G[V-\{u,v\}]$ containing the edges $v_iv_j$, $0\le i <j< \Delta$.
 We show that we can obtain a cylindrical drawing of $G_0'$ that is simple and angularly monotone in which two edges cross
if and only if this was the case in the former drawing of $G_0'$. The vertices of $G_0'$ are cyclically
 ordered around the cylinder as follows: $v_0, \ldots , v_{\Delta-1}$. An edge $v_iv_j$, $i<j$, passes in the obtained $x$-monotone drawing above
 $v_{i'}$, $i'<i$ or $i'>j$, if $v_iv_j$
crosses $vv_{i'}$ in the original drawing. Otherwise, $v_iv_j$ passes below $v_{i'}$. Moreover, the order of the intersections
of edges of $G_0'$ in the obtained drawing with the vertical line passing through $v_{i'}$ corresponds to the order of the intersections
of edges of $G_0'$ with the path $P_{i'}$  in the original drawing.

Since $G_0'$ has $\Delta$ vertices and the corresponding drawing of $G_0'$ is angularly monotone, by the hypothesis of the theorem, $G_0'$
contains a disjoint matching of size $c(\Delta)$.

 \subsection{Proof of Theorem~\ref{thm:main}}

First, we show that an angularly monotone complete simple  topological graph $G$ on $n$ vertices  contains a disjoint
matching of size at least $\Omega(\sqrt{n})$. Once this is established, by comparing $\frac n\Delta$ and $\Omega(\sqrt{\Delta})$,
and applying Theorem~\ref{thm:reduction} the claim follows.

 Consider an angularly monotone complete simple  topological graph $G$.
The vertices $v_i$ of $G$ are cyclically  ordered around the cylinder as follows: $v_0, \ldots , v_{n-1}$.

Let $G_0$ and $G_1$, respectively, denote the subgraph of $G[V]$ containing the edges $v_iv_j$, $0< i< \lfloor n/2 \rfloor <j< n$,
 such that $v_iv_j$ crosses vertical line through $v_0$ and $v_{n-1}$.
It follows that $G_0$ or $G_1$ contains  $\Omega(n^{2})$ edges. Without loss of generality, suppose that $G_0$ has more edges than $G_1$.
The graph $G_0$ can be easily transformed into an $x$-monotone simple topological graph in which two edges cross if and only if this was the case in
the cylindrical drawing of $G_0$.
Thus, by Lemma~\ref{lemma:main2}, $G_0$, and hence also $G$, contains a disjoint matching of size at least $\Omega(\sqrt{n})$.

\section{Remarks}

It is believed that the bound in  Theorem~\ref{thm:main} is not tight. This is supported by the fact that the result of G. T\'oth~\cite{T00} is believed to be
improvable to a bound linear in $kn$ and the likelihood of  Conjecture~\ref{conj:main} being true.
Therefore we find it surprising that by  methods completely different from those used by Suk~\cite{S12} we arrived at a bound with the same order of magnitude.
Another feature that both proofs have in common is a seeming difficulty of extending the proof to the case of dense graphs, i.e. the graphs with at least $\gamma n^2$ edges for some constant $\gamma>0$. A point where an adaptation of our proof for dense graphs  would fail is the application of Corollary~\ref{strong}. Indeed, Corollary~\ref{strong} is tight in the sense
that it fails if we only assume that the vertex inside the face $F$ is connected to all but one vertex incident to $F$ (regardless of the size of $F$). 
Hence, in order to extend our approach to dense graphs, a robust variant of Corollary~\ref{strong} seems to be necessary.

%

\section{Acknowledgements}

I would like to thank J\'anos Pach and Zuzana Safernov\'a for useful suggestions
concerning the presentation of the result. 

\bibliographystyle{aac}
\bibliography{refs}

\end{document}